\documentclass[12pt]{amsart}
\usepackage{xcolor}
\usepackage{graphicx}
\usepackage{geometry}

\newtheorem{theorem}{Theorem}[section]

\DeclareMathOperator{\Ric}{Ric}
\theoremstyle{definition}
\newtheorem{remark}[theorem]{Remark} 

\title[A nonisoparametric hypersurface with constant principal curvatures]{A nonisoparametric hypersurface with constant principal curvatures}
\author[A. Rodr\'{\i}guez-V\'{a}zquez]{Alberto Rodr\'{\i}guez-V\'{a}zquez}
\address{Department of Mathematics,
	University of Santiago de Compostela, Spain.}
\email{a.rodriguez@usc.es}
\thanks{The author has been supported by the project MTM2016-75897-P (AEI/FEDER, Spain) and a FPU fellowship from  Ministerio de Ciencia, Innovaci\'on y Universidades.}


\begin{document}

\begin{abstract}
In this note we construct an explicit example of a (compact) conformally flat Riemannian  manifold which admits a totally geodesic foliation of codimension one with no isoparametric leaves. This answers negatively the question: is every hypersurface with constant principal curvatures isoparametric?
\end{abstract}

\maketitle

\section{Introduction}
A hypersurface  $M$ of a Riemannian manifold $\overline{M}$ is called isoparametric if it and all its sufficiently close  (locally defined) parallel hypersurfaces  have constant mean curvature. 
Cartan proved in \cite{hipcartan}  that a hypersurface in a Riemannian manifold of constant curvature has constant principal curvatures  if and only if it is isoparametric.  In fact, there are many examples of isoparametric hypersurfaces with nonconstant principal curvatures  in spaces with nonconstant curvature \cite{damekricci}, \cite{ge}. However, it seems that there are no examples in the literature of hypersurfaces with constant principal curvatures which are not isoparametric. Even more than that, finding a nonisoparametric closed hypersurface with constant principal curvatures in a complex projective space would yield a counterexample to the longstanding Chern conjecture~\cite[Corollary~1.1]{ge}.

In this note we construct a conformally flat metric in $\mathbb{R}^n$ that admits a (non-Riemannian) foliation by totally geodesic, nonisoparametric hyperplanes. Moreover, the metric and the foliation descend to the $n$-dimensional torus $\mathbb{T}^n$. This provides an example of a nonisoparametric hypersurface with constant principal curvatures in a Riemannian manifold. Also, it shows that the equivalence between isoparametricity and constancy of the principal curvatures in spaces of constant curvature does not hold in the more general setting of conformally flat spaces.

In order to find such a metric, we need the isometry group to be sufficiently small to spoil the good behaviour of parallel hypersurfaces. Indeed, if a conformally flat space admits a transitive group of isometries, then it is locally symmetric \cite{tak}, which would lead us to the apparently outstanding problem of finding such an example in the context of symmetric spaces  \cite[\S6]{survey}. On the other hand, we construct the metric so that its isometry group is not too small so as to compute some geodesics explicitly.

\section{The ambient manifold}
Let $(x_1,\ldots,x_n)$ denote the usual coordinates in $\mathbb{R}^n$ and  $(\partial_1,\ldots,\partial_n)$ the associated coordinate vector fields.
For each  $n\geq2$ we define a metric  $$g_{ij}(x_1,\ldots,x_{n}):=h^2(x_1,\ldots,x_{n})\delta_{ij},$$
where $\delta_{ij}$ is the Kronecker's delta and $$h(x_1,\ldots,x_{n}):=\prod_{i=1}^{n-1}(2+\cos(\pi x_i))\in \mathbb{R}, \hspace{0.1cm}\mbox{for each} \hspace{0.1cm} (x_1,\ldots,x_{n})\in \mathbb{R}^n.$$
Clearly, $g$ is conformally flat. We will denote  $\mathbb{R}^n$ equipped with the metric $g$ by $\overline{M}^n$. 


%
%
\begin{remark}
In particular $g$ is invariant under translations of the lattice $2\mathbb{Z}^n$. Hence, our metric $g$ descends to the torus $\mathbb{T}^n=\mathbb{R}^n/(2\mathbb{Z}^n)$.
\end{remark}
\subsection{Christoffel symbols of $\overline{M}^n$}
It is known that Christoffel symbols are given by
\[ \Gamma_{ij}^k=\frac{1}{2}g^{kl}\big(g_{jl,i} + g_{li,j} - g_{ij,l}\big),\]
 for $i,j,k\in \{1,\ldots,n\}$, where we are using Einstein summation convention and we have denoted the partial derivative with respect to $x_i$ by $_{,i}$.
Thus,
\[\Gamma_{ij}^k=\frac{\delta_{jk}}{2h^2} h^2_{,i} + \frac{\delta_{ki}}{2h^2}h^2_{,j}-\frac{\delta_{ij}}{2h^2}h^2_{,k}. \]
Now for $n\geq2$ we have that 
\begin{subequations}\label{christoffel}
	\begin{align}
		\label{christoffel:a}
	\Gamma_{ii}^i&= (\delta_{in}-1)\frac{\pi\sin(\pi x_i)}{2+\cos(\pi x_i)}, & \Gamma_{ij}^k&=0,  \\
		\label{christoffel:b}
		\Gamma_{ij}^i&= (\delta_{jn}-1)\frac{\pi\sin(\pi x_j)}{2+\cos(\pi x_j)}, & 	\Gamma_{ii}^k&= (1-\delta_{kn})\frac{\pi\sin(\pi x_k)}{2+\cos(\pi x_k)},
	\end{align}
\end{subequations}
for mutually distinct $i,j,k\in\{1,\ldots,n\}$.
\subsection{Some vertical geodesics of  $\overline{M}^n$}

Let us define $$\Omega:=\{(a_1,\ldots,a_{n-1},x_n)\in \mathbb{R}^n: a_i\in \mathbb{Z}, \hspace{0.1cm} 0\leq i\leq n-1\}.$$ Let $a=(a_1,\ldots,a_{n-1},x_n)\in \Omega$ and $\gamma_ {a}$ be the unit-speed geodesic starting at $a$ with initial direction $\partial_n$. 
By the definition of $g$ the following maps are isometries of $\overline{M}^n$ for each $i=1,\ldots,n$
\begin{itemize}
\item $\Lambda_i: (x_{1},\ldots,x_i, \ldots, x_{n})\in \mathbb{R}^n\mapsto  (x_{1},\ldots,-x_i, \ldots, x_{n})\in \mathbb{R}^n$,
\item $\Psi_i: (x_{1},\ldots,x_i, \ldots, x_{n})\in \mathbb{R}^n\mapsto  (x_{1},\ldots,x_i+2, \ldots, x_{n})\in \mathbb{R}^n$.
\end{itemize}
Now, for each $i\in\{1,\ldots, n-1\}$, we consider the isometry $ \Psi_i^{a_i}\circ \Lambda_i$. Then, we have that  $\widetilde{\gamma}_{a}(t):= \Psi_i^{a_i}\circ \Lambda_i(\gamma_{a}(t))$ is another geodesic given by $$\widetilde{\gamma}_{a}(t)=(\gamma_{a}^1(t), \ldots, -\gamma_{a}^i(t)+2a_i  ,\ldots \gamma_{a}^n(t)).$$
But $\widetilde{\gamma}_{a}(t)$ and $\gamma_{a}(t)$ have the same initial conditions. Hence, by uniqueness we have that $\gamma_{a}^i(t)=a_i$ for each $1\leq i\leq n-1$. Observe that $h(a_1,\ldots,a_{n-1},x)=3^\rho$ for any $x\in \mathbb{R}$, where $\rho$ is the number of even entries of $(a_1,\ldots,a_{n-1})$. Thus, since $\gamma_{a}(t)$ is parametrized by arc length we get that 
\begin{equation}
\label{geodesic}
\gamma_ {a}(t)=(a_1,\ldots,a_{n-1},x_n+3^{-\rho}t).
\end{equation}

\subsection{The Jacobi operator}
\label{jacob}
It is clear that $\{\partial_i\}_{i=1}^n$ is an orthogonal global frame for $\overline{M}^n$. We will compute $\overline{R}_{\partial_n}$, the Jacobi operator associated with $\partial_n$.

All we have to do is to compute the  entries $\overline{R}_{innj}$  of the curvature tensor $\overline{R}$ for each $i,j\in\{1,\ldots,n\}$. If $i=n$ or $j=n$, then $\overline{R}_{innj}=0$. If $i,j\neq n$, then 
\begin{align*}
\overline{R}_{innj}&=\langle \overline{\nabla}_{\partial_i} \overline{\nabla}_{\partial_n}\partial_n, \partial_j\rangle - \langle \overline{\nabla}_{\partial_n} \overline{\nabla}_{\partial_i}\partial_n, \partial_j\rangle - \langle \overline{\nabla}_{[\partial_i,\partial_n]}\partial_n, \partial_j\rangle.
\end{align*}
On the one hand
\begin{align*}
\langle \overline{\nabla}_{\partial_i} \overline{\nabla}_{\partial_n}\partial_n, \partial_j\rangle&=\langle \overline{\nabla}_{\partial_i} (\Gamma_{nn}^k \partial_k), \partial_j\rangle=\langle \Gamma_{nn,i}^k{\partial_k} + \Gamma_{nn}^k\Gamma_{ik}^l \partial_l, \partial_j\rangle=h^2(\Gamma_{nn,i}^j+\Gamma_{nn}^k\Gamma_{ik}^j),
\end{align*}
and on the other hand
\begin{align*}
	\langle \overline{\nabla}_{\partial_n} \overline{\nabla}_{\partial_i}\partial_n, \partial_j\rangle&=\langle \overline{\nabla}_{\partial_n} (\Gamma_{in}^k \partial_k), \partial_j\rangle=\langle \Gamma_{in,n}^k{\partial_k} + \Gamma_{in}^k\Gamma_{nk}^l \partial_l, \partial_j\rangle=h^2(\Gamma_{in,n}^j+\Gamma_{in}^k\Gamma_{nk}^j).
\end{align*}
Since $[\partial_i, \partial_j]=0$, we conclude
\begin{equation}
\label{jacobi}
(\overline{R}_{\partial_n})_{ij}= \left \lbrace\begin{array}{lll}
h^2(\Gamma_{nn,i}^j - \Gamma_{in,n}^j +\Gamma_{nn}^k\Gamma_{ik}^j - \Gamma_{in}^k\Gamma_{nk}^j), &  \mbox{if $i,j\neq n$}\\
0, &  \mbox{in any other case.}
\end{array} \right.
\end{equation}
\section{The example}
Let $\mathcal{F}=\{\mathcal{F}_{s}\}_{s\in \mathbb{R}}$, where $\mathcal{F}_{s}=\{(x_1,\ldots,x_n)\in \mathbb{R
}^{n}: x_{n}=s\}$, for each $s\in \mathbb{R}$.
It is clear that $\mathcal{F}$ is a foliation of codimension one on $\overline{M}^n$. Let $\mathcal{S}$, $\mathcal{H}$ and $\nu \mathcal{F}_s$ denote the shape operator, the mean curvature and the normal bundle of $\mathcal{F}_s$, respectively. Then, each leaf is totally geodesic since $\partial_n\in \Gamma(\nu \mathcal{F}_s)$ and using (\ref{christoffel:a}) and (\ref{christoffel:b}), we have that  
$\langle\mathcal{S}_{\partial n} \partial_i,\partial_j\rangle=- \langle \overline{\nabla}_{\partial_i} \partial_n, \partial_j \rangle=-h^2\Gamma_{in}^j=0,$
for each $i,j=1,\ldots,n-1$.
\begin{remark}
Again, since $\mathcal{F}$ is invariant by the action of $2\mathbb{Z}^n$, this foliation descends to the torus $\mathbb{T}^n$.
\end{remark}

Let us consider $p\in \mathcal{F}_s$, $\gamma:[0,\varepsilon)\rightarrow \overline{M}^n$ a unit-speed geodesic with $\gamma(0)=p$ and $\dot{\gamma}(0)\in \nu_p \mathcal{F}_s$ for some $\varepsilon>0$, and $M^r$ the parallel hypersurface of $M$ at distance $r>0$ satisfying $\gamma(r)\in M^r$. By the Riccati equation (cf. \cite[Equation 3.8]{gray}), we have
\[\frac{d}{dr}\mathcal{S}_{\dot{\gamma}(r)}^r=\overline{R}_{\dot{\gamma}(r)} + (\mathcal{S}_{\dot{\gamma}(r)}^r)^2, \hspace{0.5cm} \mathcal{S}_{\dot{\gamma}(0)}^0=\mathcal{S}_{\partial_n},
\]
where $\mathcal{S}_{\dot{\gamma}(r)}^r$ is the shape operator of $M^r$ at $\gamma(r)$ with respect to the normal vector $\dot{\gamma}(r)$.
Now we take the trace, so
\begin{equation}
\label{traceeq}
\frac{d}{dr}\mathcal{H}_{\dot{\gamma}(r)}^r= \overline{\Ric}(\dot{\gamma}(r),\dot{\gamma}(r)) +   ||  \mathcal{S}_{\dot{\gamma}(r)}^r||^2,  \hspace{0.5cm} \mathcal{H}_{\dot{\gamma}(0)}^0=\mathcal{H},
\end{equation}
where $\mathcal{H}_{\dot{\gamma}(r)}^r$ denotes the mean curvature of $M^r$ at $\gamma(r)$, $\overline{\Ric}$ is  the Ricci tensor of $\overline{M}^n$ and $||\cdot||$  the Hilbert--Schmidt norm of a matrix.

Now we prove that no leaf of $\mathcal{F}$ is  isoparametric. Let us consider  $a\in \mathcal{F}_s\cap \Omega$ for some $s\in \mathbb{R}$.
First note that $\dot{\gamma}_a=3^{-\rho}\partial_n$  by \eqref{geodesic}. By (\ref{christoffel:a}) and (\ref{christoffel:b}), $\Gamma_{ij}^k(\gamma_{a}(t))=0$ and $\Gamma_{in}^n=0$ . Hence, by (\ref{jacobi}),  we have
\[\overline{\Ric}(\dot{\gamma}_a(r),\dot{\gamma}_a(r)) =\sum_{i=1}^{n-1} \Gamma_{nn,i}^i(\gamma_{a}(t))=\pi^2(1-n+\frac{4}{3}\rho), \]
where $\rho$ is the number of even entries of $(a_1,\ldots,a_{n-1})$.

As a consequence,  if $a=(0,\ldots,0,s)\in \mathcal{F}_s\cap \Omega$ and $b=(1,\ldots,1,s)\in \mathcal{F}_s\cap \Omega$, we get that
\begin{align*}
\overline{\Ric}(\dot{\gamma}_a(r),\dot{\gamma}_a(r))=\frac{n-1}{3}\pi^2>0& &\mbox{and} & & \overline{\Ric}(\dot{\gamma}_b(r),\dot{\gamma}_b(r))=(1-n)\pi^2<0.
\end{align*}

But in our case, $|| \mathcal{S}_{\dot{\gamma}_{a}(0)}||^2=|| \mathcal{S}_{\dot{\gamma}_{b}(0)}||^2=0$. Therefore, by (\ref{traceeq}), we deduce that $\frac{d}{dr}\rvert_{ r=0}\mathcal{H}_{\dot{\gamma}_{a}(r)}^r>0$ and $\frac{d}{dr}\rvert_{ r=0}\mathcal{H}_{\dot{\gamma}_{b}(r)}^r<0$. This way we can conclude that, for  small $r>0$ the mean curvature of the parallel hypersurface of $\mathcal{F}_s$ at distance $r>0$,  is not constant. Then, $\mathcal{F}_s$ is not isoparametric.


\begin{thebibliography}{99}
	\bibitem{hipcartan} \'E. Cartan: Familles de surfaces isoparam\'etriques dans les espaces \`a courbure constante, \textit{Ann. Mat. Pura Appl.} \textbf{17} (1938), no. 1, 177–191.
	\bibitem{damekricci}	J. C. D\'iaz-Ramos, M. Dom\'inguez-V\'azquez:
	 Isoparametric hypersurfaces in Damek--Ricci spaces, \textit{Adv. Math.} \textbf{239} (2013), 1--17. 
	\bibitem{survey} J. C. D\'iaz-Ramos, M. Dom\'inguez-V\'azquez, V. Sanmart\'in-L\'opez: Submanifold geometry in symmetric spaces of noncompact type, \textit{ S\~ao Paulo J. Math. Sci.}, (2019), doi:10.1007/s40863-019-00119-6, 1-36. 
	\bibitem{ge}J. Ge, Z. Tang,  W. Yan: A filtration for isoparametric hypersurfaces in Riemannian manifolds, \textit{J. Math. Soc. Japan}, \textbf{67} (2015), no. 3, 1179–1212.
	\bibitem{gray} A. Gray: \textit{Tubes},   Progress in Mathematics \textbf{221}, Birkh\"auser Basel, Boston, (2004).
	\bibitem{tak}H. Takagi: Conformally flat Riemannian manifolds admitting a transitive group of isometries II, \textit{T\^ohoku Math. J.}, \textbf{27} (1975), 445-451. 
	
\end{thebibliography}
\end{document}